\documentstyle{article}
 \title{Fano Manifolds, Contact Structures,
and Quaternionic Geometry}
\author{Claude LeBrun\thanks{Supported in part
by NSF Grant DMS-9204093.}\\
Department of Mathematics\\ SUNY, Stony Brook, NY 11794-3651}
\def\be{\begin{equation}}
\def\ee{\end{equation}}
\def\bea{\begin{eqnarray*}}
\def\eea{\end{eqnarray*}}

\def\hook{\mbox{}\begin{picture}(10,10)\put(1,0){\line(1,0){7}}
  \put(8,0){\line(0,1){7}}\end{picture}\mbox{}}
\newcounter{exam}[section]
\renewcommand{\theexam}{\thesection.\arabic{exam}}
\newenvironment{xpl}{\medskip \refstepcounter{exam}
\noindent {\bf Example \theexam .}}{\hfill $\diamondsuit$\mbox{}\bigskip}
\newcounter{remark}[section]
\renewcommand{\theremark}{\thesection.\arabic{remark}}
\newenvironment{rmk}{\medskip \refstepcounter{remark}
\noindent {\bf Remark \theremark .}}{\hfill $\Box$\mbox{}\bigskip}
\newtheorem{main}{Theorem}
\newtheorem{thm}{Theorem}
\newtheorem{defn}{Definition}

\newtheorem{prop}{Proposition}
\newtheorem{cor}{Corollary}

\newtheorem{maincor}[main]{Corollary}

\newcommand{\sctn}{\setcounter{equation}{0}
\setcounter{thm}{0} \setcounter{prop}{0}
\setcounter{defn}{0}\setcounter{cor}{0}
\setcounter{exam}{0} \setcounter{remark}{0}
\setcounter{lem}{0}\section}
\newenvironment{proof}{\medskip
\noindent {\bf Proof.}}{\hfill \rule{.5em}{1em}\mbox{}\bigskip}
\def\Bbb{\bf}

\def\P{{\Bbb CP}}
\def\C{{\Bbb C}}
\newcommand{\eel}[1]{\label{#1}\end{equation}}

\def\bp{{\Bbb P}}

\begin{document}
\maketitle
\begin{abstract} Let $Z$ be a compact complex
$(2n+1)$-manifold which carries a {\em complex contact structure},
meaning a codimension-1 holomorphic sub-bundle
$D\subset TZ$ which is maximally non-integrable.
If $Z$ admits a K\"ahler-Einstein metric
of positive scalar curvature, we show that  it is
the Salamon twistor space of a quaternion-K\"ahler
manifold $(M^{4n}, g)$. If $Z$
also admits a second complex contact structure
$\tilde{D}\neq D$, then $Z= {\Bbb CP}_{2n+1}$.
 As an application, we give
several new characterizations of the
Riemannian manifold ${\Bbb HP}_n=Sp(n+1)/\left(
Sp(n)\times Sp(1)\right)$.
\end{abstract}

 \sctn{Introduction}
If $(M,g)$ is an oriented Riemannian $m$-manifold, how many tensor
fields  $\varphi\neq 0$ does $M$ admit which
satisfy  $\nabla \varphi =0$,
where  $\nabla$ is the
Levi-Civit\`a connection associated with the
 Riemannian metric $g$? Obviously there are always
some such fields:  the metric,
the volume form, their inverses, and linear combinations of
contracted tensor products
of these. There will be more
only if the geometry of
$M$ is special in the sense that the {\em holonomy group} --- that
is \cite{bes,Sbook}, the group of linear transformation of
a tangent space induced by parallel transport around loops --- is
a proper subgroup of $SO(m)$, and a
 fundamental research topic for
present-day Riemannian geometry is the problem of
classifying those  Riemannian manifolds
 whose holonomy geometries are special in this sense.
  In fact, if we exclude
the locally symmetric spaces (for which the curvature
tensor ${\cal R}$ satifies $\nabla {\cal R}=0$) and local products
of manifolds of lower dimension, there are (up to conjugation)
 only seven possible
families of  connected
holonomy groups:
$SO(m)$, $U(m/2)$, $SU(m/2)$,
$G_2$ ($m=7$), $\mbox{Spin}(7)$ ($m=8$),
 $Sp(m/4)$, and
$[Sp(m/4)\times Sp(1)]/{\Bbb Z}_2$ ($m\geq 8$).
The  present paper is largely motivated by the
problem of understanding the last of these holonomy families.

\begin{defn} A Riemannian manifold $(M, g)$
of dimension $4n$, $n\geq 2$, is said to be
quaternion-K\"ahler if its holonomy group
is (conjugate to) a  subgroup of $[Sp(n)\times Sp(1)]/{\Bbb Z}_2
\subset SO(4n)$. \label{qk}
\end{defn}
This terminology may be justified by the fact \cite{bes,Sbook} that
any such manifold carries a 4-form   $\varphi$ such
that  $\nabla \varphi =0$,
and this 4-form may be thought of as a quaternionic
analogue of the
 2-form which plays such a central r\^ole in the K\"ahler
geometry of complex manifolds.

A quaternion-K\"ahler manifold is necessarily Einstein \cite{ber},
and its scalar curvature is non-zero iff the
holonomy group contains the $Sp(1)$ factor of
  $[Sp(n)\times Sp(1)]/{\Bbb Z}_2$. In particular,
such a manifold is necessarily compact provided it is complete and
its (constant) scalar
curvature is positive.

\begin{defn}
A quaternion-K\"ahler manifold $(M,g)$
will be called {\em positive} if it is
compact and has positive scalar curvature.
\end{defn}
Of course, multiplying the metric $g$ of such a manifold
by a positive constant always gives us  `new'
examples in a completely trivial way.
We henceforth  eliminate this  possibility by
normalizing the scalar curvature
of $g$ to be the same as that of
the standard metric on ${\Bbb HP}_n=S^{4n+3}/Sp(1)$,
namely $s=16n(n+2)$.

It should be emphasized  that a quaternion-K\"ahler
manifold
is typically  {\em not} a K\"ahler manifold;
$[Sp(n)\times Sp(1)]/{\Bbb Z}_2\not \subset U(2n)$.
Nonetheless, questions about such
manifolds can, in the positive case,
 be reduced to problems
in algebraic geometry over $\Bbb C$ by means of
a   {\em twistor construction}
 \cite{sal,beber} which associates
a  compact complex $(2n+1)$-manifold $Z$  with any
positive quaternion-K\"ahler $4n$-manifold $(M,g)$.
Moreover, the so-called
twistor spaces $Z$ that arise by this construction
have some remarkable  properties.
First, they admit
K\"ahler-Einstein metrics of
positive scalar curvature;
 in particular, they have $c_1>0$, and so are
{\em Fano manifolds}.
Secondly, any such $Z$ admits a {\em complex contact
structure}, meaning a maximally non-integrable
holomorphic sub-bundle $D\subset TZ$ of
complex codimension 1.

In short, any twistor space $Z$
is a {\em Fano contact manifold} in the sense of  Definition
\ref{fan} below. On the other hand, every known
 Fano contact manifold is
actually a twistor space. Our first
main result indicates that this is no accident:

\begin{main} \label{links}
Let $Z$ be a Fano contact manifold.
Then $Z$ is a twistor  space iff it
admits a K\"ahler-Einstein metric.
\end{main}
Based on this, it might seem reasonable to conjecture that
every Fano contact manifold is a twistor space.
For other evidence in favor of such a conjecture, see
\cite{lfin,lsal,ps,wolf}.

 The same machinery used to prove the above
also allows us to address
 the  problem of determining  which
Fano contact manifolds carry more than one contact
structure:

\begin{main}\label{only}
 Let $Z$ be a Fano   manifold of complex dimension
$2n+1$ which admits a K\"ahler-Einstein metric.
  If $Z$ admits two distinct complex contact structure
$D, \tilde{D}\subset TZ$,
then $Z\cong{\Bbb CP}_{2n+1}$.
\end{main}
This result has some  interesting ramifications
 regarding quaternion-K\"ahler manifolds.
To give a typical
 example, let us say that  a self-diffeomorphism
$\Phi : M\to M$ of a quaternion-K\"ahler manifold $(M,g)$
is a {\em quaternionic automorphism} if it preserves the
$[GL(n,{\Bbb H})\times Sp(1)]/{\Bbb Z}_2$ structure
determined by $g$; thus the
 projective  transformations of ${\Bbb HP}_n$
induced by
 elements of $GL(n+1, {\Bbb H})$
are quaternionic automorphism of
${\Bbb HP}_n$,
whereas only   elements of the subgroup
$Sp(n+1)\times {\Bbb R}^+\subset
GL(n+1, {\Bbb H})$ act by isometries.
However,  Theorem \ref{only} implies that
this example is essentially anomalous:

\begin{maincor} Let $(M,g)$ be a positive quaternion-K\"ahler
$4n$-manifold.
If there is a quaternionic automorphism $\Phi: M\to M$ which does not
preserve $g$,
 then $(M,g)$ is isometric to the symmetric space ${\Bbb HP}_n$.
\end{maincor}

\vfill

\noindent {\bf Notation.} Certain
notational conventions  employed in these pages
may elicit gasps from the occasional
 algebraic geometer who takes the time to peruse them.
Thus multiplicative rather than additive notation will
be employed for the Picard group, and  our convention
regarding projectivizations of vector spaces and vector
bundles is that  ${\Bbb P}(E)=(E-0)/{\Bbb C}^{\times}$.
On the other hand, both the real and holomorphic tangent bundles
of a complex manifold $Z$ will generally be denoted by
$TZ$ except in cases  where this might be
 likely to cause confusion.

\bigskip

\noindent {\bf Acknowledgements.}
The  author would like to express his gratitude to
 the Newton Institute of Cambridge University
for its  hospitality and financial support.
He would also like to thank   Charles Boyer,
Krzysztof Galicki, and
Yun-Gang Ye for the discussions
which provided the
original impetus for this research.

\pagebreak

\sctn{Complex Contact Manifolds} \label{cont}

\begin{defn} A complex contact manifold is a pair $(Z,D)$, where
$Z$ is a complex manifold  and
$D\subset TZ=T^{1,0}Z$ is a codimension-1 holomorphic sub-bundle
which is maximally non-integrable in the sense that the
 O'Neill tensor (`Frobenius obstruction')
\begin{eqnarray*}
D\times D &\to& TZ/D\\
(v,w)&\mapsto & [v,w]\bmod D
\end{eqnarray*}
is everywhere  non-degenerate.
\end{defn}

The condition of non-integrability has a very useful reformulation,
which we shall now describe. Given a codimension-1 holomorphic sub-bundle
$D\subset TZ$, let $L:=TZ/D$ denote the quotient line bundle;
we thus have an exact sequence
\be 0\to D\longrightarrow TZ\stackrel{\theta}{\longrightarrow} L\to
0\label{one}\ee
where $\theta$ is the
tautological  projection. But
 we may also think of $\theta$ as a line-bundle-valued
1-form
$$\theta \in \Gamma (Z, \Omega^1(L))~,$$
and so attempt to form its exterior derivative $d\theta$.
Unfortunately, this ostensibly depends on a choice of local
trivialization; for if $\vartheta$ is any 1-form,
$d(f\vartheta)=fd\vartheta +df\wedge\vartheta$. However, it is now clear that
$d\theta|_D$ {\em is} well defined as a section of $L\otimes \wedge^2D^{\ast}$,
and an elementary computation, which we leave to the reader, shows that
 $d\theta|_D$, thought of in this way, is exactly the  O'Neill tensor mentioned
above. Now if the skew form $d\theta|_D$ is to be non-degenerate, $D$ must
have positive even rank $2n$, so that  $Z$ must have odd complex dimension
$2n+1\geq 3$.
Moreover, the non-degeneracy exactly requires that
\be \theta \wedge (d\theta )^{\wedge n}\neq 0 \label{ndg}
\ee
as a section of $
\Omega^{2n+1}(L^{n+1})$.
This now provides a bundle isomorphism between
$L^{\otimes (n+1)}$
and the anti-canonical line bundle $K^{-1}=\wedge^{2n+1}T^{1,0}Z$,
in keeping with the isomorphism
$D\cong L\otimes D^*$ induced by the O'Neill tensor.

Conversely, let $Z$ be a   compact complex $(2n+1)$-manifold
with $H^1(Z, {\Bbb Z}_{n+1})=0$,
and suppose that $c_1(Z)$ is divisible by $n+1$. Then there is a
unique holomorphic line bundle $L:=K^{-1/(n+1)}$ such that
$L^{\otimes (n+1)}\cong K^{-1}$. If we are then given a
twisted holomorphic 1-form
$$\theta \in \Gamma (Z, \Omega^1(K^{-1/(n+1)}))$$
we may then construct
$$\theta \wedge (d\theta )^{\wedge n}\in
\Gamma (Z, \Omega^{2n+1}(K^{-1}))=
\Gamma (Z,{\cal O})= {\Bbb C}~.$$
If this constant is non-zero, one calls
$\theta$  a {\em complex contact form}
because the corresponding
 $D=\ker \theta$ is then a complex
contact structure. Now the holomorphic map
$$ \Gamma (Z, \Omega^1 (L))\ni \theta
\mapsto\theta\wedge (d\theta )^n \in{\Bbb C} $$
is homogeneous of degree $(n+1)$, and so
is represented by a homogeneous polynomial.
Since two elements of $\Gamma (Z, \Omega^1 (L))$
define the same sub-bundle $D$ iff they
are  proportional, it follows that the space of contact
structures, if non-empty,
 is the complement of
a degree-$(n+1)$ complex hypersurface
$S\subset {\Bbb P}[\Gamma (\Omega^1 (L))]$;
in particular, {\em the space of all
contact structures on $Z$ is a connected complex
manifold}.

\begin{xpl} \label{proj}
Let ${\Bbb V}\cong \C^{2n+2}$ be an
even-dimensional  complex
vector space, and let
$\Upsilon\in \wedge^2{\Bbb V}^*$ be any skew form on
${\Bbb V}$. Then there is an associated
$\theta \in \Gamma ({\Bbb P}({\Bbb V}), \Omega^1(2))$
defined by $\varpi^*\theta:= \Xi\hook \Upsilon$, where
$\varpi : {\Bbb V}-0\to {\Bbb P}({\Bbb V})$ is the canonical
projection, and where $\Xi$ is the (Euler) vector
field which generates the ${\Bbb C}^{\times}$-action of
scalar multiplication on $\Bbb V$; conversely, every
element of $\Gamma ({\Bbb P}({\Bbb V}), \Omega^1(2))$ arises
this way from a $\Upsilon$, as  may, for example,
be read off
from the Euler exact sequence
$$ 0\to \Omega^1(2) \to {\Bbb V}^*(1)\to {\cal O}(2)\to 0$$
on ${\Bbb P}({\Bbb V})$. On the other hand,
the anti-canonical line bundle $K^{-1}$ of
${\Bbb P}({\Bbb V})$ is isomorphic to  ${\cal O}(2n+2)$,
so we have $K^{-1/(n+1)}\cong {\cal O}(2)$, and
any contact structure on ${\Bbb P}({\Bbb V})\cong {\Bbb CP}_{2n+1}$
must arise from some $\theta \in \Gamma ({\Bbb P}({\Bbb V}), \Omega^1(2))$.
Now the condition that $\theta\wedge (d\theta)^n\neq 0$
may now be rewritten as $\Upsilon^{n+1}\neq 0$
as a consequence of the observation that
$d(\varpi^*\theta)=d (\Xi\hook \Upsilon)= \pounds_{\Xi}
{\Upsilon}=2\Upsilon$. Thus the set of contact structures
on  ${\Bbb CP}_{2n+1}$
is exactly parameterized by the symplectic
forms on $\C^{2n+2}$ modulo rescalings.
\end{xpl}

\begin{xpl} \label{cot}
 Let $Y_{n+1}$ be any complex manifold, and let
$Z_{2n+1}={\Bbb P}(T^{\ast}Y)$
be its projectived
 holomorphic cotangent bundle. Observe that
the cotangent bundle $T^*Y$ carries a
tautological 1-form $\hat{\theta}$ defined by
$\hat{\theta}|_{\vartheta}= p^*\vartheta$,
where $p: T^*Y\to Y$ is the
canonical projection. If $\Phi_t: T^*Y\to T^*Y$
is fiber-wise  scalar multiplication
by $t\in {\Bbb C}^{\times}$,
we also have $\Phi_t^*\hat{\theta}=t\hat{\theta}$,
so there is  a line-bundle-valued
form $\theta\in \Gamma ({\Bbb P}(T^*Y), \Omega^1(L))$
on $Z$ such that $\varpi^*\theta=\hat{\theta}$; here
$\varpi$ is the canonical projection from $T^*Y-0_Y$
to $Z={\Bbb P}(T^{\ast}Y)$, and $L\to Z$
is the holomorphic line bundle whose local sections
are homogeneity-1 functions on $T^*Y-0_Y$.
On the other hand, $\Upsilon=d{\hat{\theta}}$
is a holomorphic symplectic form on
$T^*Y$, and the non-degeneracy of $\Upsilon$
implies that $\theta$ is a contact form.
In particular, $L^{n+1}$ is isomorphic to the
anti-canonical line bundle of $Z= {\Bbb P}(T^{\ast}Y)$.
\end{xpl}

The usual proof of the Darboux theorem \cite{arn} for real contact
manifolds applies equally well in the complex
case; thus,  any complex contact manifold $(Z,D)$ of
dimension $2n+1$  is locally
isomorphic to
$({\Bbb C}^{2n+1}, \ker (dz^{2n+1}+\sum_{j=1}^n z^jdz^{j+n}))$.
As a consequence, any complex contact manifold
may be obtained by gluing together open sets in
${\Bbb C}^{2n+1}$ with transitions functions
which are {\em complex contact transformations} ---
i.e. biholomorphisms which preserve the
fixed complex contact structure.

In order to better
understand this notion of complex contact
transformation, we should first try to understand
the infinitesimal version; that is,
when does the pseudo-group of
biholomorphisms generated by a holomorphic vector field
preserve a given contact structure $D$?
The answer is given by the following result:

\begin{prop}
Let $(Z,D)$ be any complex contact manifold,
and let ${\bf u}\in \Gamma (Z, {\cal O}(L))$ be any
holomorphic section of the contact line bundle.
Then there is a unique holomorphic vector
field $\zeta$ on $Z$ such that
$\theta (\zeta ) = {\bf u}$ and such that
the pseudo-group of local transformations of
$Z$ generated by $\zeta$ consists of
complex contact transformations.\label{contra}
\end{prop}
\begin{proof}
It suffices to prove the result locally, since the
uniqueness will guarantee that local choices
of $\zeta$ agree on overlaps; thus
 we may assume that the
contact line bundle $L$ is trivial.
Relative to a  trivialization
of $L$, the contact form $\theta$ becomes
a holomorphic 1-form
$\vartheta$,
whereas the section ${\bf u}$ is represented by an ordinary
holomorphic function $u$. Now there is a unique holomorphic
vector field $\eta$ such that $\eta\hook d\vartheta =0$
and $\eta \hook \vartheta = 1$ because
 $\vartheta\wedge (d\vartheta)^n\neq 0$ and
the condition that $\theta (\zeta ) = {\bf u}$
now becomes $\zeta= u\eta + \xi$ for some
$\xi \in \Gamma ({\cal O} (D))$.
On the other hand, the
condition that $\zeta$ generate
contact transformations can be written as
$\pounds_{\zeta}\vartheta= f\vartheta$
for some holomorphic function $f$,
and we therefore must have
$$ (u\eta + \xi) \hook d\vartheta  + d [( u\eta + \xi) \hook \vartheta ]
\equiv 0 \bmod \vartheta ~ , $$
and hence
$\xi \hook d\vartheta = -du \bmod \vartheta$.
Thus
\be \zeta = u \eta - (d\vartheta |_D)^{-1} (du |_D)\eel{expl}
is the unique infinitesimal contact transformation with
$\theta (\zeta ) = {\bf u}$.
\end{proof}

\begin{cor} \label{gpsp}
The exact sequence
\be 0\to {\cal O}(D) \to {\cal O}(TZ)\stackrel{\theta}{\to }
{\cal O} (L)\to 0\eel{basic}
splits as a sequence of sheaves of  {\em abelian groups}.
In particular,
$$H^p(Z, {\cal O}(TZ))\cong
H^p (Z, {\cal O}(D))\oplus H^p(Z, {\cal O}(L))$$
for all $p$.
\end{cor}

\begin{rmk} The canonical  splitting
${\cal O}(TZ)\leftarrow {\cal O}(L)$
of (\ref{basic})
given by Proposition \ref{contra} is  {\em not}
a homomorphism  of $\cal O$-modules, since
(\ref{expl}) manifestly
involves the first derivative of ${\bf u}$.
In particular, Corollary \ref{gpsp} does
{\em not} assert the existence of a splitting
of the associated exact
sequence (\ref{one}) of vector bundles.
In fact, we will see in Corollary \ref{doesnt} that the
latter  never splits if
$Z$ is Fano.
\end{rmk}

\begin{prop} Let $(Z,D)$ be a compact complex
contact $(2n+1)$-manifold with
$H^1(Z, {\Bbb Z}_{n+1})=0$.
 Then $Z$ admits a second
contact structure $\tilde{D}\neq D$ iff
$\Gamma (Z, {\cal O}(D))\neq 0$. \label{gamd}
\end{prop}
\begin{proof}
If $\Gamma (Z, {\cal O}(D))\neq 0$,
there is
a non-trivial holomorphic vector field  $\zeta$
on $Z$ with $\theta (\zeta )=0$. Since we also have $\theta (0)=0$,
Proposition \ref{contra} asserts that  this $\zeta$ is not
an infinitesimal contact transformation. Thus the contact
structure
$(\exp t\zeta )_*D$ will differ from $D$ if $t$
is sufficiently  small.

Conversely, if there is a second contact structure
$\tilde{D}\neq D$  on $Z$, we must have $\tilde{D} = \ker
\tilde{\theta}$, $\tilde{\theta}\in \Gamma  (Z, \Omega^1 (L))$,
 by virtue of the assumption
 that $H^1(Z, {\Bbb Z}_{n+1})=0$; and since
$\tilde{D}\neq D$, there must be a point of
$Z$ at which the values of $\theta$ and $\tilde{\theta}$
are linearly independent. Thus  $\tilde{\theta}|_D\in \Gamma (Z, {\cal
O}(D^*\otimes L))$
is not  zero, and this of course implies
$(d\theta |_D)^{-1}(\tilde{\theta}|_D)\in \Gamma (Z, {\cal O}(D))$
is not zero, either.
\end{proof}

Because complex contact structures are locally trivial, and
infinitesimal contact automorphisms are specified  by
holomorphic  sections of the contact line bundle
$L$, the space of infinitesimal
 deformations of a compact complex
manifold $(Z,D)$ is exactly  $H^1(Z, {\cal O}(L))$.
Thus the existence of the natural splitting
$H^1(Z, {\cal O}(TZ))\leftarrow  H^1(Z, {\cal O}(L))$
of Corollary \ref{gpsp}  may be interpreted
as saying that
{\sl  you can't deform the contact structure without
deforming the complex structure}. The following
result \cite{lsal,nit} is therefore an immediate
consequence of the connectedness of the space of
contact forms:

\begin{prop} Let $Z$ be a simply connected
compact complex   manifold.
Then any two complex contact structures on $Z$ are
equivalent via some biholomorphism of $Z$.
\end{prop}

\begin{rmk}
There is no real analogue of this; for example,
the 3-sphere $S^3$ carries \cite{ben} many inequivalent
real contact structures. The essential  difference
 is that the set of real contact forms
is typically  disconnected.
\end{rmk}

Examples \ref{proj} and \ref{cot} were
both displayed in a way that
related them  to symplectic structures
on ${\Bbb C}^{\times}$-bundles over the
given manifold $Z$. As you might expect, this
is a manifestation of a general procedure
 \cite[Appendix 4E]{arn}
known as  {\em symplectification}, which
we now
review. Let $L^{*\times}$ denote the complement
 of the zero section in
in the total space of $L^*=L^{-1}$,
and let $\varpi : L^{*\times} \to Z$ denote the
canonical projection, which we will
consider to be a holomorphic principal
${\Bbb C}^{\times}$-bundle. Then
$\varpi^*L$ is canonically trivial, although
the canonical non-zero section has homogeneity 1
with respect to the ${\Bbb C}^{\times}$-action.
It follows that $\varpi^*\theta$ may be considered
to be a holomorphic 1-form on $L^{*\times}$ rather
than just a section  of $\Omega^1(\varpi^*L)$.
We may therefore define $\Upsilon=d(\varpi^*\theta)$,
and observe that the non-degeneracy (\ref{ndg})
of $\theta$ is exactly equivalent to requiring that the
closed holomorphic 2-form
$\Upsilon $ satisfy $\Upsilon^{n+1}\neq 0$;
that is, $\Upsilon $ is a {\em holomorphic symplectic
form}, and in particular the map $v\to \Upsilon (v, \cdot )$
is  an  isomorphism $TL^{*\times}\to T^*L^{*\times}$.
Conversely, any symplectic form of homogeneity
1 on a principal ${\Bbb C}^{\times}$-bundle arises from
a complex contact structure on the base by the
formula $\varpi^*\theta = \Upsilon ( \xi , \cdot )$,
where $\xi$ is the vector field which
generates the ${\Bbb C}^{\times}$-action.

This has interesting consequences back down on the contact
 manifold $Z$.

\begin{prop}
Let $(Z,D)$ be a complex contact manifold, and
let $L=TZ/D$ be its contact line bundle.
Then the contact form $\theta\in \Gamma (Z, \Omega^1(L))$
gives rise to a non-degenerate section $\Upsilon\in
\Gamma (Z, L^*\otimes \wedge^2 J^1L)$,
where $J^1L$ is the 1-jet bundle of
$L$. In particular, there is an isomorphism
$(J^1L)^*\to L^*\otimes J^1L$ induced by contraction with
$\Upsilon$. Moreover, this is compatible with
the  map  $D\otimes L^*\to  D^*$
induced by $d\theta|_D$. \label{symp}
\end{prop}
\begin{proof} Since a local holomorphic section $f$ of
$L$ may be identified with  a function $\hat{f}$ on $L^{*\times}$
which has homogeneity 1 with respect to the
${\Bbb C}^{\times}$-action, the identification
$J^1(f)\leftrightarrow d\hat{f}$ gives us a
${\Bbb C}^{\times}$-invariant identification of
$\varpi^*L^*\otimes J^1L$ with the cotangent bundle of $L^{*\times}$.
Thus the symplectic form $\Upsilon=d(\varpi^*\theta)$ on
$L^{*\times}$,
which transforms under the ${\Bbb C}^{\times}$-action
with homogeneity 1, may be identified with a
non-degenerate holomorphic section of
$L\otimes \wedge^2(L^*\otimes J^1L)=L^*\otimes \wedge^2J^1L$.
Moreover, the diagram
$$\begin{array}{ccc}
J^1L& \stackrel{\Upsilon^{-1}}{\longrightarrow} & L\otimes (J^1L)^*
 \\ \uparrow & & \downarrow\\
\Omega^1(L ) & \longrightarrow &TZ\\
\downarrow & & \uparrow\\
D^*\otimes L&\stackrel{(d\theta |_D)^{-1}}{\longrightarrow} & D
\end{array}$$
commutes as a consequence of the fact that
$\Upsilon = d(\varpi^*\theta )$.
\end{proof}

\begin{rmk}  Proposition \ref{contra} says that
$\Gamma (Z, {\cal O}(L))$ can be identified
with the set of infinitesimal contact transformations
of $(Z,D)$, and so has a Lie algebra structure.
In fact, the Lie bracket $\Gamma ({\cal O}(L))\times
\Gamma ({\cal O}(L))\to \Gamma ({\cal O}(L))$ is
exactly given by $[{\bf u}, {\bf v}]=\Upsilon^{-1}(J^1{\bf u},
J^1{\bf v})$. In the same way, one can also give
$\bigoplus_{m=1}^{\infty}\Gamma (Z, {\cal O}(L^m))$
the structure of a graded Lie algebra.
\end{rmk}

\begin{prop} Let $Z$ be a compact complex $(2n+1)$-manifold
of K\"ahler type,
and let $D$ be a compact complex structure on $Z$.
Then the obstruction $\in {\bf Ext}^1_Z({\cal O}(L), {\cal O}(D))=
H^1(Z, {\cal O}(L^*\otimes D))$ to splitting
the  exact sequence
$$0\to D \to TZ\to L\to 0$$
 is obtained  by applying
the composition $$H^1(\Omega^1)\to H^1( {\cal O}(D^*))
\stackrel{(d\theta|_D)^{-1}}{\longrightarrow}
 H^1( {\cal O}(L^*\otimes D))$$
to $\frac{2\pi i}{n+1}
c_1(Z)\in H^1(Z, \Omega^1)$.
\label{extclass}
\end{prop}

\begin{proof}
The obstruction to splitting the jet sequence
$$ 0\to \Omega^1(L)\to J^1L \to L\to 0 $$
is \cite{at} the Atiyah obstruction $a(L)\in H^1(Z, \Omega^1)$,
and may be expressed in \v{C}ech cohomology
as $[d\log f_{\alpha\beta}]$, where
$\{ f_{\alpha\beta}\}$ is a system of transition functions
for $L$. Thus the image of $a(L)$ in
$H^1(Z,  {\cal O}(L^*\otimes D))$ is the
extension class of
\be 0\to D^*\otimes L \to (J^1L)/{\cal O}\to L\to 0 .\eel{jet}
But contraction with $\Upsilon^{-1}$
 converts (\ref{jet}) into the
exact sequence
 \be 0\to D  \to TZ\to L \to 0 .\label{two} \ee
The result thus follows from
Proposition \ref{symp} and the observation \cite{at}
that $a(L)= 2\pi i c_1(L)$ for any line bundle
on a compact K\"ahler manifold.\end{proof}

\begin{cor} \label{doesnt}
Suppose that $(Z,D)$ is a compact complex
contact manifold such that $c_1(Z)>0$. Then (\ref{one})
does not split.
\end{cor}
\begin{proof}
Because $K\otimes L$ is a negative
line bundle, the Kodaira vanishing
theorem implies that $H^1(Z, {\cal O}(L^*))=0$.
Thus the restriction map
$H^1(Z, \Omega^1)\to H^1(Z, {\cal O}(D^*))$
is injective, and the result
follows from Proposition \ref{extclass}.
\end{proof}

While a number of interesting
things can be said
regarding  complex contact manifolds in general,
 one seems to need strong
extra hypotheses before a classification becomes
imaginable. One possible route is to
limit ones ambitions to the low dimensional cases; cf.
\cite{ye2}. In this article, however, we will instead
study manifolds of the following type:

\begin{defn} If $(Z,D)$ is a complex contact manifold
such that $c_1(Z)>0$,
we will say that $(Z,D)$ is a
{\em Fano contact manifold}. \label{fan}
\end{defn}

Thus, for example, ${\Bbb CP}_{2n+1}$ and
${\Bbb P}(T^*{\Bbb CP}_{n})$ proxide two examples
of Fano contact manifolds fo dimension $2n+1$.
These are both examples of {\em homogeneous}
complex contact manifolds, meaning that their
groups of contact transformations act transitively;
there is \cite{wolf} exactly one such object for each
simple complex Lie algebra,
and every  homogeneous complex contact
manifold is   automatically Fano
as a consequence of Proposition \ref{contra}.
 From our perspective, however,  the overwhelming
reason to study such objects is provided by the
following result \cite{sal,beber}:

\begin{thm}[Salamon/B\'erard-Bergery]
Let $(M^{4n},g)$ be a
quaternion-K\"ahler manifold of scalar curvature
$16n(n+2)$. Then there is
a complex contact $2n+1$-manifold $(Z,D)$,
called the {\em twistor space} of $(M,g)$,  which
admits a K\"ahler-Einstein metric $h$ of scalar curvature
$8(n+1)(2n+1)$
such that
\begin{description}\item{(i)} there is a  Riemannian submersion
 $\wp : Z\to M$ with totally geodesic fibers
  $S^2$ of constant curvature $4$;
\item{(ii)} the horizontal sub-bundle of this submersion is
the contact distribution $D\subset TZ$;
\item{(iii)} each fiber of $\wp$ is a  rational complex  curve
${\Bbb CP}_1\subset Z$, with
normal bundle   isomorphic to $[{\cal O}(1)]^{\oplus 2n}$;
and
\item{(iv)} there is a free anti-holomorphic involution
$\sigma : Z\to Z$ which commutes with $\wp$.
\end{description}
\label{sal}
\end{thm}

In order to allow for the case $n=1$ while ensuring
that this theorem remains true,
we now supplement Definition \ref{qk} as follows:

\begin{defn} An oriented Riemannian manifold $(M, g)$
of dimension $4$ is said to be
quaternion-K\"ahler if it is
Einstein and has self-dual Weyl curvature. \end{defn}

Now one might hope that this twistor-theoretic
machinery would provide a powerful source of
new examples of Fano contact manifolds.
However, the only known positive quaternion-K\"ahler
manifolds are symmetric spaces, and their
twistor spaces are exactly the homogeneous complex
contact manifolds alluded to above. Indeed, the
following result \cite {lfin,lsal}
gives one reason to wonder whether there
are non-homogeeous examples at all:

\begin{thm}
There are, up to biholomorphism, only finitely
many Fano contact manifolds $(Z, D)$ of any
fixed dimension.
\end{thm}

There are also results \cite{hit,ps}
that show that there are no
non-homogeneous examples with $n=1,2$.
And in all dimensions one has the following result
\cite{lfin,lsal}:

\begin{thm} Let $Z$ be a  contact Fano  manifold.
If $b_2(Z)\geq 2$, then $Z\cong \bp(T^*{\P_{n+1}})$.
\label{first}
\end{thm}

\pagebreak

\sctn{Contact Structures and Einstein Metrics}

The following well-known
observation is apparently due to
Calabi \cite{cal}.

\begin{prop}
Let $(Z,h)$ be a K\"ahler-Einstein manifold of
positive scalar curvature and complex dimension $m$;
and  let $K^{\times}\to Z$ be the
${\Bbb C}^{\times}$-bundle obtained from the canonical
line bundle $K$ by deleting the zero section.
Then $K^{\times}$ carries an incomplete Ricci-flat K\"ahler metric
${\bf h}$ such that $\Phi_t^*{\bf h}= |t|^{2/(m+1)}{\bf h}$,
where $\Phi_t: K\to K$
is scalar multplication by $t\in {\bf C}^{\times}$.
Moreover, $(Z,h)$ is a K\"ahler quotient of $(K^{\times},{\bf h})$
by
 $S^1\subset {\Bbb C}^{\times}$.
\label{rf}
\end{prop}

\begin{proof}
Define $r: K\to K$  by $r(x)= \frac{1}{2}\| x\|^{2/m+1}$, where
$\|\cdot\|$ is the norm on $K$ associated with the
K\"ahler-Einstein metric and $m$ is the complex dimension of
$Z$. We may then define ${\bf h}$ to be the K\"ahler
metric on $K^{\times}$ associated with
the positive $(1,1)$-form $\omega_1:=i\partial \bar{\partial} r$.

We may suppose that the K\"ahler-Einstein metric $h$ is mormalized
so that its Ricci tensor is $2(m+1) h$. The
 pull-back of $\omega$ via the canonical
projection  $\varpi: K^{\times}\to Z$ is then
$$\varpi^*\omega = \frac{i}{2}\partial\bar{\partial} \log r , $$
so that  $(Z,\omega)$ is precisely
 the symplectic quotient of $(K^{\times},\omega_1)$ corresponding to
$r=1/2$.
On the other hand, there is a tautological holomorphic
$(m,0)$-form $\phi$ on $K^{\times}$
defined by $\phi|_{\varphi}=\varpi^*\varphi$, and
the relation between  $r$ and
the norm $\|\cdot\|$ on $K\to Z$  tells us that
$$\varpi^*\omega^m\propto
\frac{\phi\wedge \bar{\phi}}{r^{m+1}},$$
where $\propto$ means the two expressions differ
by multiplication by a non-zero constant. On the other
hand, if $\Xi$ denotes the holomorphic (Euler)
vector field which generates the ${\Bbb C}^{\times}$-action
on $K^{\times}$, then $\Xi\hook \partial \phi = \pounds_{\Xi}\phi=
\phi$, whereas $\Xi\hook (\partial\log r \wedge \phi)
= (\Xi \log r) \phi = \phi /(m+1)$; thus
$\partial \phi = (m+1) \partial \log r \wedge \phi$,
since both sides have type  $(m+1,0)$
and $m+1=\dim_{\Bbb C} K^{\times}$. Hence
\bea (\omega_1)^{m+1}&=& (i\partial\bar{\partial} r)^{m+1}\\
&=& \left(2r\varpi^*\omega + i\frac{\partial r
\wedge \bar{\partial}r}{r}\right)^{m+1}\\
&\propto&  r^{m-1} \varpi^* \omega^m\wedge
\partial r\wedge \bar{\partial r}\\
&\propto&r^{m-1}\frac{\phi\wedge \bar{\phi}}{r^{m+1}}
\wedge \partial r\wedge \bar{\partial} r\\
&=& \phi\wedge \bar{\phi}
\wedge \partial \log r\wedge \bar{\partial} \log r\\
&\propto& \partial \phi\wedge\overline{\partial \phi}.
\eea
Since $\partial\phi$ is a holomorphic form,
the logarithm of the volume form of $\omega_1$ is thus
pluri-harmonic, and the K\"ahler metric $\bf h$ is Ricci-flat.
\end{proof}

Descending back to $Z$, we thus have the following:

\begin{prop} Let $Z$  be a Fano manifold which
  admits a K\"ahler-Einstein metric.
Then the 1-jet bundle $J^1K^*\to Z$ of   the
anti-canonical line  bundle of $Z$ admits a natural Hermitian-Einstein
inner product induced by the Ricci-flat metric $\bf h$
on $K^{\times}$. In particular, $J^1K^*$ is  quasi-stable ---
that is, it is a semi-stable direct sum of stable vector bundles.
\label{stab}
 \end{prop}
\begin{proof}
First observe that the pull-back of  $K\otimes J^1K^*$
from $Z$ to $K^{\times}$
is ${\Bbb C}^{\times}$-equivariantly isomorphic to
$T^*K^{\times}$,
since any
local holomorphic  section $f$ of $K^*$ on $Z$ can
be identified with a  holomorphic function $\hat{f}$
on $K^{\times}$ of
homogeneity $1$, and the value of the 1-jet of $f$
at a point of $Z$ determines and is determined by
the value of $d\hat{f}$ at any point of the
corresponding fiber of $K^{\times}$. In other words,
a local section of $J^1K^*$ on $Z$ is a
holomorphic cotangent field $\psi$ on $K^{\times}$ which
satisfies $\Phi_t^*\psi = t\psi$ for all
$t\in {\Bbb C}^{\times}$. Since
the inner product $\langle \cdot , \cdot \rangle$
on covectors determined by
$\bf h$ satisfies
 $\Phi_t^*\langle \cdot , \cdot \rangle =
 |t|^{-2/(m+1)}\langle \cdot , \cdot \rangle$, we
can now define an inner product $(\cdot , \cdot )$
on $J^1K^*\to Z$ by
$$(\psi , \tilde{\psi}) :=
r^{-2m/(m+1)} \langle \psi , \tilde{\psi}\rangle .$$
We claim that this inner product is Hermitian-Einstein.

Since this claim is local in character, we may now restrict
to an open set in $Z$ over which we have
a  root $\ell=K^{-1/(m+1)}$ of the
the anti-canonical bundle. Since the
curvature of the obvious Chern connection
on this root is a constant multiple
of the K\"ahler form $\omega$, it
is thus sufficient for us to check that
the induced inner product on
$\ell^{-m}\otimes J^1K^*\cong J^1\ell$
is Hermitian-Einstein. However, a local
section of this bundle may be interpreted
as a local holomorphic cotangent vector
field $\psi$  on $K^{\times}$ such that $\Phi_t^*\psi= t^{1/(m+1)}\psi$,
and the induced inner product on such objects
is just $\langle \cdot , \cdot \rangle$. Thus
the curvature of $(T^*K^{\times}, \langle \cdot , \cdot \rangle)$
is just the pull-back of the curvature
of $\ell^{-m}\otimes J^1K^*$; in particular,
$T^*K^{\times}$ is flat along the fibers of $K^{\times}\to Z$. On the
other hand, $\varpi^*\omega \equiv
 \omega_1/2r \bmod (\partial r , \bar{\partial} r)$,
so   the projection
$\varpi : K^{\times}\to Z$ is conformal in the horizontal directions.
Combining these last two observations tells us
that the Ricci curvature $F^{\alpha}_{\beta j\bar{k}}\omega^{j\bar{k}}$
of $\ell^{-m}\otimes J^1K^*$ pulls back to a multiple
of the Ricci curvature of $\bf h$, and so vanishes.
Hence the inner product induced on $J^1K^*$ by
$\bf h$ is Hermitian-Einstein, as claimed.
The quasi-stability statement now follows from the work of
Kobayashi  and L\"ubke \cite{kob,lueb}.
\end{proof}

\begin{rmk}
The extension class of
$$0\to {\cal O} \to   (K\otimes J^1K^*)^*\to TZ\to 0$$
is just  $c_1(Z)\in H^1(Z, \Omega^1)=Ext^1_Z(TZ, {\cal O})$,
so the stability aspect of the above result
should be attributed to Tian  \cite{tian}. However, see
 Remark \ref{error} below.
\end{rmk}

This stability result now allows us to prove the following:

\begin{thm} Let $(Z,D)$ be a compact complex contact manifold,
and suppose that  $h$ is
  a K\"ahler-Einstein metric
 of positive scalar curvature on $Z$.
Then the Ricci-flat manifold $(K^{\times},{\bf h})$ of
Proposition \ref{rf} is finitely covered by a hyper-K\"ahler
manifold.
\label{hypK}
\end{thm}

\begin{proof}
Let $m=2n+1$ be the complex dimension of $Z$, and let
$L=K^{-1/(n+1)}$ be the contact line bundle.
Since $J^1K^*=L^n\otimes J^1L$ is Hermitian-Einstein by
Proposition \ref{rf}, it follows that $L^*\otimes \wedge^2 J^1L=
L\otimes K^2\otimes \wedge^2 J^1K^*$ is also Hermitian-Einstein.
But the symplectification procedure of page \pageref{symp}
gives us a non-degenerate section $\Upsilon$ of
$L^*\otimes\wedge^2 J^1L\subset {\cal H}om ( (J^1L)^*, L^*\otimes J^1L)$,
so  contraction with
$\Upsilon^{-1}\otimes \Upsilon^{-1}$ is an
isomorphism
$$L^*\otimes \wedge^2J^1L\stackrel{\cong}{\longrightarrow}
(L^*\otimes \wedge^2J^1L)^* .$$
Thus the odd Chern classes of $L^*\otimes \wedge^2 J^1L$
are all 2-torsion, and in particular  this bundle
has degree 0.
This shows that
 $L^*\otimes \wedge^2 J^1L$
is actually Hermitian-Ricci-flat, meaning that
its curvature $\hat{F}$ satisfies $\omega\cdot \hat{F}=0$.

Given a contact form $\theta\in\Gamma (Z, \Omega^1(L))$,
the symplectification
procedure just alluded to produces  a non-degenerate
$\Upsilon \in \Gamma (Z,{\cal O}(L^{-1}\otimes \wedge^2J^1L))$, and, since
$L^{-1}\otimes \wedge^2J^1L$ is Hermitian-Ricci-flat,
  $\Upsilon $ must be parallel.
Indeed, if $\nabla$ is the Chern connection
of $L^{-1}\otimes \wedge^2J^1L$, we have
$\nabla^{0,1}\Upsilon =0$ because $\Upsilon $ is
holomorphic, and whereas
\bea\int_Z\|{\nabla}^{1,0} \Upsilon \|^2~d\mu&=&
 -\int_Z \langle \Upsilon ,
i\omega^{j\bar{k}}\nabla_{\bar{k}}\nabla_j\Upsilon \rangle~d\mu\\&=&
\int_Z \langle \Upsilon ,
i\omega^{j\bar{k}}\hat{F}_{j\bar{k}}(\Upsilon ) \rangle~d\mu\\&=&0.
\eea
 Hence $\nabla \Upsilon =
\nabla^{1,0}\Upsilon +\nabla^{0,1}\Upsilon =0$.

Let $L^{*\times}$ be the complement of the
zero section in $L^*$, and notice that there is a
natural covering map $L^{*\times}\to K^{\times}$  given by
 $x\to x^{\otimes (m+1)}$, and recall that
the symplectification of procedure of page \pageref{symp}
initially displays
$\Upsilon $ as a non-degenerate holomorphic
2-form on $L^{*\times}$, in a manner consistent with
our identification of sections of
$J^1L$ with homogeneity-1 1-forms on $K^{\times}$.
Thus,  pulling $\bf h$ back to
$L^{*\times}$ via this covering,
$\Upsilon $ becomes a parallel $(2,0)$-form on  the
Ricci-flat K\"ahler manifold $(L^{*\times}, {\bf h})$.
The holonomy group $\subset SU(2n+2)$ of $\bf h$ thus
stabilizes a $(2,0)$-form of maximal rank, and so
must be a subgroup of $Sp(n+1)$.
Hence $(L^{*\times}, {\bf h})$ is a  hyper-K\"ahler manifold.
\end{proof}

\begin{rmk}
A more direct way of seeing  that
$L^*\otimes \wedge^2J^1L$ is Hermitian-Ricci-flat
is to observe that, in terms of the
local root $\ell$ used in the
proof of Proposition \ref{stab},
one has $L^*\otimes \wedge^2J^1L=\wedge^2 J^1\ell$.
It's curvature is thus explicitly given
by ${\hat{F}^{\alpha \beta}}_{\gamma \epsilon j\bar{k}}=
2\delta^{ [ \alpha}_{ [ \gamma }
F^{ \beta ]}_{\epsilon ] j\bar{k}}$,
and so is annihilated by contraction
with $\omega^{j\bar{k}}$.
\end{rmk}

\begin{thm}
 Let $(Z,h)$ be a compact K\"ahler-Einstein manifold
with positive scalar curvature. If $Z$
 admits  two distinct contact
structures $D$ and $\tilde{D}$, then the
associated Ricci-flat manifold
$(K^{\times}, {\bf h})$ is actually flat. \label{flat}
\end{thm}

\begin{proof}
The same argument used in
the proof of Theorem \ref{hypK} says that the holonomy
group of $(L^{*\times}, {\bf h})$ now stabilizes {\em two}
linearly independent $(2,0)$-forms of maximal rank;
moreover, these two $(2,0)$-forms are homogeneous
of degree 1 with respect to the ${\Bbb R}^+$-action
generated by $\xi=\mbox{grad\/}r$.
Thus $L^{*\times}$ can then  be covered by open sets $U$
which split as  Riemannian products
$(U_1, {\bf h}_1)\times (U_2, {\bf h}_2)$ of pairs
hyper-K\"ahler manifolds; moreover, these local splittings
 fit together to yield an orthogonal pair of
 foliations of $L^{*\times}$, and these foliations are
${\Bbb R}^+$-invariant. But the
   homogeneity of  ${\bf h}$ says that
$\pounds_{\xi } {\bf h} = 2 {\bf h}$,
and  its components $\xi_1$ and $\xi_2$
tangent to
$U_1$ and $U_2$ therefore satisfy
$\pounds_{\xi_1}{\bf h}_1= 2{\bf h}_1$
and $\pounds_{\xi_2}{\bf h}_2= 2{\bf h}_2$; that is,
 $\xi$, $\xi_1$, and $\xi_2$ are
{\em homothetic vector fields}. If we now let
$S\subset L^{*\times}$ denote the real hypersurface $S=1$,
then $\xi$  is a unit normal vector field on $S$,
and the homothetic nature of $\xi$ implies that the
second fundamental form of $S$ is just the restriction $\bf g$
of $\bf h$ to $S$. If now we defined $f: S\to {\bf R}$ to be
$|\xi_1|^2$,  then $f$ is not identically zero and
$|df|=|\xi_1| ~ |\xi_2|$;
hence $\xi_1$ and $\xi_2$ respectively have zeroes
at the minima and maxima of $f$. Now if
$p\in S$ is a point at which at which $\xi_1=0$,
let us consider the diffeomorphism $\exp \xi_1 : U_1'\to U_1$
induced on a small leaf-wise
neighborhood $U_1'\subset U_1$ of $p$. Since
$(\exp \xi_1)^*{\bf h}_1= e^2{\bf h}_1$,
its derivative $(\exp \xi_1)_{*p}: T_pU_1\to T_pU_1$
at $p$ is diagonalizable over $\Bbb C$, with all
 eigenvalues of modulus $e^2> 1$.
Hence all the eigenvalues of the
induced push-forward maps on
$T_pU_1\otimes \bigotimes^{k+3} (T_p^*U_1)$
have modulus $e^{-2(k+2)} < 1$. However,
since $(\exp \xi_1)_*$ just multiplies ${\bf h}_1$ by a
constant, it preserves the covariant derivatives
$\nabla^{(k)}{\cal R}:= \nabla \cdots \nabla {\cal R}$
of the curvature tensor of ${\bf h}_1$; and since $1$ is not an
eigenvalue of
$$(\exp \xi_1)_{*p}:T_pU_1\otimes {\bigotimes}^{k+3} (T_p^*U_1)\to
T_pU_1\otimes {\bigotimes}^{k+3} (T_p^*U_1),$$ it follows that
 $\nabla^{(k)}{\cal R}$
vanishes  at $p$ for all $k$. Thus all the $U_1$ components
of the curvature tensor of $\bf h$ on $U=U_1\times U_2$
vanish to infinite order at $p$, and since $\bf h$ is
real-analytic, it follows that all its
first-foliation curvature components
are identically zero on  $L^{*\times}$.
Since the same argument can also be applied to $\xi_2$,
 we conclude that ${\bf h}$ is actually flat.
\end{proof}

\setcounter{main}{1}
\begin{main} Let $Z$ be a compact complex
 $(2n+1)$-manifold  which admits
a K\"ahler-Einstein metric of positive scalar curvature.
Then $Z$ admits two distinct  complex contact structures
iff $Z\cong {\Bbb CP}_{2n+1}$.
\end{main}

\begin{proof}
By Theorem \ref{hypK},
the Ricci-flat K\"ahler metric
$\bf h$ is
flat, whereas the proof also tells us that the
real hypersurface $S$ given by $r=1/2$ has
second fundamental form equal to the induced
metric $\bf g$. By the Gauss-Codazzi equations,
$(S, {\bf g})$  therefore has constant curvature
$1$, and so is locally isometric to the
standard unit sphere in ${\Bbb R}^{4n+4}$.
If $\eta$ denotes the vector field which
generates the circle action $x\to e^{it}x$
on $S$, then $\eta$ is a unit Killing vector
field on $S$, and so corresponds to
a $(4n + 4) \times (4n + 4)$ matrix which is
simultaneously in $so(4n+4)$ and $SO(4n+4)$.
Such a matrix is an isometric complex
structure on ${\Bbb R}^{4n+4}$.
We may therefore isometrically identify
any sufficiently small open set in $S$
 with an open set in the unit sphere in
 $S^{4n+3}\subset{\Bbb C}^{2n+2}$
in such a way that $\eta$ corresponds to the generator of
multiplication by $e^{it}$.
But the fibration $S\to Z$ is a Riemannian
submersion, and the fibers are just the orbits
of the $U(1)$-action generated by $\eta$. Hence
$Z=S/U(1)$ is locally isometric to
the Fubini-Study metric on ${\Bbb CP}_{2n+1}=S^{4n+3}/U(1)$.
The developing map thus gives us an open isometric  immersion
from  the universal cover of $(Z,h)$ to the symmetric space
${\Bbb CP}_{2n+1}$; and since $Z$ and ${\Bbb CP}_{2n+1}$
are both compact and simply connected,
this developing map therefore
defines a global isometry between $(Z,h)$
and the symmetric space  ${\Bbb CP}_{2n+1}$.
Since there is only one complex structure
which is compatible with the Fubini-Study metric,
this isometry is  necessarily a
biholomorphism, and we are done.
\end{proof}

\begin{rmk}
By essentially the same argument, one can also sharpen
Proposition \ref{stab} as follows: {\em
if $(Z,h)$ is a compact K\"ahler-Einstein $m$-manifold
of positive scalar curvature, and if $J^1K_Z^*$ is
not strictly stable, then $Z\cong {\Bbb CP}_m$. }
\end{rmk}

With Proposition \ref{gamd}, Theorem \ref{only} now implies

\begin{cor} Suppose that  $Z_{2n+1}\not \cong {\Bbb CP}_{2n+1}$
is a Fano manifold with complex contact structure $D$.
Also assume that  $Z$ admits a K\"ahler-Einstein metric $h$.
Then $\Gamma(Z, {\cal O}(D))=0$.
\label{prelim}
\end{cor}

\begin{cor} Let  $(Z,D)$ be a Fano contact manifold
which admits a K\"ahler-Enstein metric. Then
the following are equivalent:
\begin{description}
\item{(a)} $Z\cong{\Bbb CP}_{2n+1}$;
\item{(b)}  $h^0(Z, \Omega^1 (L)) > 1$;
\item{(c)} $Z$ admits more than one complex contact structure;
\item{(d)} $Z$ has an automorphism which is not a contact transformation;
\item{(e)} $h^0(Z, {\cal O}(D))>0$;
\item{(f)} $\chi (Z, {\cal O} (D))> 0$.
\end{description}
\end{cor}

\begin{proof} For  any Fano contact manifold $(Z,D)$ one has
$H^p(Z, {\cal O}(D)) = 0$  $\forall p > 1$ as a consequence
\cite{lsal,sal}
of the
Kodaira vanishing theorem. Thus, invoking Theorem
\ref{only} and the proof of Proposition \ref{contra},
  (f) $\Rightarrow$ (e) $\Rightarrow$ (d) $\Rightarrow$ (c)
$\Rightarrow$
(a) $\Rightarrow$ (f).  On the other hand, the equivalence of
(b) and (c) is obvious.
\end{proof}

\begin{rmk} Let $Z$ be a Fano manifold of
complex dimension $m$ and Picard number $b_2(Z)=1$.
Tian \cite[Theorem 2.1]{tian} claims to prove that if
the tangent bundle of  $Z$ admits
a proper holomorphic sub-bundle $D\subset TZ$,
then  \be\frac{c_1(D)}{c_1(Z)}
 \leq \frac{\mbox{rank} D}{m+1},\eel{ineq}
 with
equality only if there is a non-trivial
holomorphic section of $D\subset TZ$. However, the complex
contact structure $D$ of any twistor space $Z_{2n+1}$
 is  a sub-bundle with
$\frac{c_1(D)}{c_1(Z)}= \frac{n}{n+1}=
\frac{\mbox{rank} D}{m+1}$. Thus,
by  Corollary \ref{prelim},
 the twistor space of
$\widetilde{Gr}_4({\Bbb R}^{n+4})=SO(n+4)/SO(n)\times SO(4)$
provides a counter-example to this assertion,
as does the twistor space of any  exceptional Wolf space \cite{wolf}.
Fortunately, the mistake seems to occur  only in
the last line of the proof, where it is implicitly assumed that
the contraction of two non-zero tensors is necessarily
non-zero. The proof of
(\ref{ineq}) is
thus unaffected.\label{error}
\end{rmk}

\setcounter{main}{0}
\begin{main}
Let $Z$ be a Fano contact manifold.
Then $Z$ is the twistor  space of a
quaternion-K\"ahler manifold iff it
admits a K\"ahler-Einstein metric.
\end{main}
\begin{proof} Suppose that $Z$ is K\"ahler-Einstein
and admits a complex contact structure.
By Theorem \ref{hypK},
the incomplete Ricci-flat metric $\bf h$ on
$L^{*\times}$ admits a non-degenerate parallel $(2,0)$-form
$\Upsilon $ of homogeneity 1, and in particular
is hyper-K\"ahler. Expressing such a
form in terms of its real and
imaginary parts, we have $\Upsilon =\omega_3+i\omega_2$, and
using the metric $\bf h$ transform the
2-forms $\omega_1$, $\omega_2$, into skew endomorphisms
${\bf J}_2$ and ${\bf J}_3$ of $TL^{*\times}$, we have
${\bf J}_1{\bf J}_2=-{\bf J}_2{\bf J}_1={\bf J}_3$
and ${\bf J}_3{\bf J}_1=-{\bf J}_1{\bf J}_3={\bf J}_2$,
where ${\bf J}_1$ denotes the tautological complex structure of
$L^{*\times}$. But since ${\bf J}_2\neq 0$ is a parallel infinitesimal
orthogonal
transformation of $TL^{*\times}$, ${\bf J}_2^2$ has a negative eigenvalue,
and by replacing $\Upsilon $ with a suitable real multiple, we
may assume that this eigenvalue is $-1$. The
$(-1)$-eigenspaces of ${\bf J}_2^2$ now form a parallel sub-bundle
of $L^{*\times}$ which is invariant under the ${\Bbb R}^+$-action
generated by $\xi$. If this is a proper-sub-bundle,
the proofs of Theorems \ref{flat} and
\ref{only} then  show that $Z$ is the twistor space
${\Bbb CP}_{2n+1}$ of ${\Bbb HP}_n$. Otherwise,
${\bf J}_2^2=-1$ and hence ${\bf J}_3^2=-1$, so that ${\bf J}_A$,
 $A=1,2,3$, is an anti-commuting triple of compatible
complex structures for our hyper-K\"ahler metric $\bf h$,
and the corresponding K\"ahler forms
$\omega_A$ satisfy $\pounds_{\xi} \omega_A= 2\omega_A$,
  $A=1,2,3$, as a consequence of the homogeneity of
$\Upsilon $ ($A=2,3$) and our definition of
$\bf h$ ($A=1$).
But since the 2-forms  $\omega_A$ are closed,
this may be rewritten as
$d(\xi\hook\omega_A)=2\omega_A$; and as
 $\xi\hook\omega_A= {\bf h}({\bf J}_A\mbox{grad\/} r, \cdot)
={\bf J}_Ad r$, we thus have
$$i\partial_A\bar{\partial}_Ar=
{\textstyle \frac{1}{2}} d\/ {\bf J}_Adr=\omega_A, ~A=1,2,3.$$
  But these three equations
exactly say that $r$ is a {\em hyper-K\"ahler potential}
in the sense of Swann \cite{swann}, and \cite[Proposition 5.5]{swann}
now tells us that the ${\bf J}_A\xi$ are normalized generators of an
isometric $sp(1)$-action on the hypersurface $S\subset L^{*\times}$
defined by $r=1/2$. Since $S$ is compact, this action exponentiates
to yield $Sp(1)$-action,  and since,  by construction,
 the normalized generator
${\bf J}_1\xi$ has period $\pi$, this
action descends to an action of  $SO(3)=Sp(1)/{\Bbb Z}_2$.
And since, by construction, the $SO(2)$-action
generated by ${\bf J}_1\xi$ is free,
this $SO(3)$-action is free, as a consequence
of the transitivity of  the adjoint action
of  $SO(3)$ on the unit sphere in $so(3)$.
 Thus $M:=S/SO(3)=L^{*\times}/{\Bbb H}^{\times}$ is a smooth
compact $4n$-manifold,
and the Riemannian submersion metric $g$ on $M$ is
\cite[Theorem 5.1]{swann}
quaternion-K\"ahler, with twistor space  $Z$.
\end{proof}

\begin{rmk}
Instead of invoking Swann's results \cite{swann}
on the relation between
hyper-K\"ahler and quaternion-K\"ahler manifolds, we could have
instead proceeded by showing that the triple
of vector fields ${\bf J}_A\xi$ give the hypersurface
$S$ a so-called {\em  3-Sasakian structure}
\cite{ish,kon}. Long relegated
to obscurity and neglect, this concept
 has recently turned out to be
 a remarkably fruitful
source  of compact  Einstein manifolds \cite{bgm}.
 Theorem \ref{only} is similarly related to \cite{tayu}.
\end{rmk}

\begin{rmk} Uwe Semmelmann has brought some related
results \cite{spin1,spin2} to our
attention which,
when $Z$ is K\"ahler-Einstein and  {\em spin},
link the presence of a
 complex contact structure to that of a
{\sl Killing spinor}; when $Z$ is spin, this would
provide yet another strategy for proving Theorem \ref{links}.
Note, however, that
 a  Fano contact $2n+1$-manifold
$Z\neq {\Bbb CP}_{2n+1}$
is spin iff $n$ is odd; thus this method
would only prove
 `half' of the result in question.
\end{rmk}

The following result is well-known; cf.  \cite{BE,bes,L0,nit}.

\begin{cor} Suppose that the
twistor space $Z$ of the
 positive quaternion-K\"ahler
manifold $(M,g)$ is biholomorphic to ${\Bbb CP}_{2n+1}$.
Then $(M,g)\cong {\Bbb HP}_n$.
\label{pin}
\end{cor}

\begin{proof}
Since $h^0(Z, \Omega^1(K^{-1/(n+1)}))=
h^0({\Bbb CP}_{2n+1}, \Omega^1(2))=(n+1)(2n+1)>1$,
$Z$ has more than one complex contact structure.
The  Swann hyper-K\"ahler metric $\bf h$ associated with
Salamon's K\"ahler-Einstein metric $h$ on $Z$
is therefore flat by Theorem \ref{flat}. The totally
umbilic hypersurface given by $r=1/2$ is therefore
locally isometric to the unit sphere
in ${\Bbb R}^{4n+4}$, and the three
unit Killing fields ${\bf J}_A\xi$
define three orthogonal complex
structures on this ${\Bbb R}^{4n+4}$
which generate a representation of
$sp(1)$. Hence $S$ is locally
isometric to $S^{4n+3}\subset {\Bbb H}^{n+1}$
in a manner which sends these three
Killing fields to right multiplication by
$i$, $j$, and $k$. Thus $M$ is locally
isometric to ${\Bbb HP}_n=S^{4n+3}/Sp(1)$, and
the developing map construction therefore
produces  a
global isometry $M\to {\Bbb HP}_n$.
\end{proof}

 Theorem \ref{only}
now yields
a much better proof of \cite[Theorem 3.2]{lsal}:

\begin{cor} Let $(M,g)$ and
$(\tilde{M}, \tilde{g})$ be two
positive quaternion-K\"ahler manifolds,
 and  let $Z$ and $\tilde{Z}$  be
their respective twistor spaces.
Then $(M,g)\cong (\tilde{M}, \tilde{g})$ as
Riemannian manifolds iff $Z \cong \tilde{Z}$
as complex manifolds.
\end{cor}

\begin{proof}  By Corollary \ref{pin},
we may
assume henceforth that $Z,\tilde{Z}\not\cong  {\Bbb CP}_{2n+1}$.
 Theorem \ref{only} thus tells us that
the only complex contact structures on
$Z$ and $\tilde{Z}$, respectively, are
 the horizontal distributions of the
twistor projections $\wp: Z\to M$ and
 $\tilde{\wp}:\tilde{Z}\to \tilde{M}$.

Now suppose that $Z$ and $\tilde{Z}$ are biholomorphic.
 The Bando-Mabuchi uniqueness
 theorem for K\"ahler-Einstein metrics \cite{BM}
then implies that that there is a biholomorphic
{\sl isometry} $\Phi$ between the
the K\"ahler-Einstein
manifolds $(Z,h)$ and $(Z,\tilde{h})$.
But because $Z$ and $\tilde{Z}$
have only one complex contact structure apiece,
the isometry $\Phi$ sends horizontal
subspaces of $\wp$ to horizontal subspaces of $\tilde{\wp}$,
and hence fibers of $\wp$ to fibers of $\tilde{\wp}$.
Thus $\Phi:Z\to \tilde{Z}$ covers a diffeomorphism
$F: M\to \tilde{M}$.
But since the twistor projections are
Riemannian submersions, $F$ is automatically
an isometry, and $(M,g)\cong (M,\tilde{g})$.

The converse follows from
 the naturality of the Salamon construction.
\end{proof}

Our next application concerns the following:

\begin{defn}\label{qaut}
Let $(M,g)$ be a quaternion-K\"ahler $4n$-manifold,
and suppose that $F : M\to M$ is a diffeomorphism.
We will say that $F$ is a {\em quaternionic
automorphism} if  the derivative of $F$
preserves the $[GL(n, {\Bbb H})\times Sp(1)]/{\Bbb Z}_2$
determined by the quaternion-K\"ahler structure.
\end{defn}

\setcounter{main}{2}
\begin{maincor}
 If $(M,g)\not\cong {\Bbb HP}_n$ is
a  positive quaternion-K\"ahler
$4n$-manifold, any quaternionic automorphism of
$M$ is an isometry.  \end{maincor}

\begin{proof}  Since the construction of the complex
manifold $Z$ only  involves \cite{sal2,BE}
 the $GL(n, {\Bbb H}) Sp(1)$-structure
of $M$, any quaternionic automorphism $F:M\to M$ lifts uniquely
to a biholomorphism $\Phi : Z\to Z$, and this biholomorphism will
 be contact iff   \cite{nit,L0} $F$ is an isometry.
But if $\Phi$ is {\em not} contact, $\Phi_*D$ is a
second contact structure, and the result follows from
Theorem \ref{only} and Corollary \ref{pin}.
\end{proof}

Another interesting consequence concerns the spectrum of the Laplacian:

\begin{cor} Let $(M,g)$ be a compact quaternion-K\"ahler
$4n$-manifold of   scalar curvature $s$.
Let $\lambda_1$ denote the
smallest positive eigenvalue of the Laplace operator
$\Delta : C^{\infty}(M)
\to C^{\infty}(M)$. Then
$$\lambda_1 \geq \frac{(n+1)s}{2n(n+2)} ~ ,$$
with equality iff $M$ is  ${\Bbb HP}_n$
and $g$ is a multiple of the standard metric.
\end{cor}

\begin{proof} Since the statement is empty if $s\leq 0$,
we may assume that $s>0$; and by
scale invariance,
we may therefore assume  henceforth that $g$ has scalar curvature
$16n(n+2)$.  In conjunction with our usual convention
that the K\"ahler-Einstein  metric $h$
on $Z$ has Ricci tensor $4(n+1)h$, this
 makes the twistor projection
 $\wp: Z\to M$ into a Riemannian
submersion with totally geodesic fibers.
If $f$ is an eigenfunction of $\Delta$ on $M$, with
eigenvalue $\lambda > 0$, $\hat{f}=\wp^*f$
is an eigenfunction of $\Delta$ on $Z$, also  with  the
eigenvalue $\lambda$; moreover, $\mbox{grad\/}\hat{f}$ is
perpendicular to the fibers of $\wp$, and so is
a section of $D$.  But
$\Delta(\Delta-8(n+1))=4\nabla_{\nu}\nabla^{\bar{\mu}}
\nabla_{\bar{\mu}}\nabla^{\nu}$
on $Z$, and the latter elliptic
operator is manifestly non-negative;
since $\hat{f}$ is
an eigenfunction of this operator, with eigenvalue
$\lambda(\lambda - 8(n+1))$, it follows that
$\lambda\geq 8(n+1)$. Moreover \cite{mats}, $\hat{f}$
is  in the kernel of this operator only if
$\mbox{grad\/}\hat{f}$ is a real-holomorphic vector field,
so that $\lambda = 8(n+1)$ only if $D$ admits a non-trivial
holomorphic section. Thus, by Corollary \ref{prelim},
$\lambda > 8(n+1)$ unless $M\cong {\Bbb HP}_n$.
Conversely, $8(n+1)$ actually
{\em is} in the spectrum of ${\Bbb HP}_n$,
since all of the many holomorphic sections of $D\subset
T{\Bbb CP}_{2n+1}$ are \cite{mats} of the form
$\mbox{grad\/}\hat{f}_1+J\mbox{grad\/}\hat{f}_2$, where
$\hat{f}_1$ and $\hat{f}_2$ are
of eigenfunctions of $\Delta$
with eigenvalue
$8(n+1)$. \end{proof}

\end{document}